\documentclass{amsart}
\usepackage{amsmath}
  \usepackage{paralist}
  \usepackage{graphics} 
  \usepackage{epsfig} 
 \usepackage[colorlinks=true]{hyperref}
\hypersetup{urlcolor=blue, citecolor=red}

\newcommand{\imagepath}{}

\title[Detection of time-varying heat sources]{Detection of time-varying heat sources using an analytic forward model}

\author{Janne P.~Tamminen}


 \email{janne.tamminen@imaqen.com}

\thanks{This work was supported by Suomen Kulttuurirahasto via the Post Docs in Companies (PoDoCo) program and done partly while working at University  of Helsinki, department of mathematics and statistics.}

\newcommand{\abs}[1]{\left\vert {#1} \right\vert}

\newcommand{\norm}[1]{\left\Vert {#1} \right\Vert}

\def\Im{\mathop{\rm Im}}
\def\Re{\mathop{\rm Re}}

\newcommand{\ad}{\"a}

\newcommand{\R}{{\mathbb R}}

\newcommand{\rmi}{\mathrm{i}}
\newcommand{\dee}{\mathbf{d}}

\date{\today}

\begin{document}

\maketitle
\tableofcontents

\begin{abstract}
We present a simple, analytic point source model for both static and time-varying point-like heat sources and the resulting temperature profile that solves the heat equation in dimension three. Simple algorithms to detect the location and spectral content of these sources are developed and numerically tested using Finite Element Mesh simulations. The resulting framework for heat source reconstruction problems, which are ill-posed inverse problems, seems promising and warrants for future research. Possible fields of application for our work are material testing, to detect manufacturing defects, and medical imaging to detect abnormal health conditions.
\end{abstract}

\section{Introduction}
This work stems from the field of medical infrared (IR) imaging, where the boundary temperature of an object is measured using an IR camera and subsequently analysed using e.g. statistical methods or thermodynamic modeling. More specifically, the starting point of this paper is the method of dynamic infrared imaging (DIRI), where time-varying temperature phenomena are recorded with framerates of 30Hz and higher. Thus we are dealing with a time-varying temperature distribution on the boundary of the object, from where we are interested in computing some thermal properties inside the object e.g. heat conductivity, blood perfusion or the heat source distribution. In this paper we deal with the last case in 3D: given the temperature measurement $T(x,t)|_{\partial \Omega}$ on the boundary of the domain $\Omega\subset\R^3$ in a given time interval, reconstruct the heat source 
\begin{equation}\label{heatSource}
q(x,t) = s(t)\delta(x_0) \quad \textrm{ in }\Omega,
\end{equation}
where $\delta(x_0)$ refers to a point source (delta-distribution) at the spatial location $x_0$ and $s(t)$ is a signal of the form
\begin{equation}\label{signal}
s(t) = \sum_{n=0}^N S_n\cos(2\pi f_n t + \phi_n).
\end{equation}
We specifically investigate two problems:
\begin{enumerate}
\item Given the static measurement $T(x)|_{\partial \Omega}$, reconstruct the location $x_0$ and magnitude $Q$ of a static heat source $q(x) = Q\delta(x_0)$.
\item Given the dynamic temperature measurement $T(x,t)|_{\partial \Omega}$, reconstruct the location $x_0$, the amplitudes $S_n$ (and the phases $\phi_n$ should we need them).
\end{enumerate}
Mathematically speaking, reconstructing the heat source from boundary measurements is an ill-posed inverse problem, which is very sensitive to noise in the measurement and thus needs a regularized solution method. Inverse problem algorithms can be categorized, roughly speaking, in three: analytic or direct methods, iterative optimization methods, and statistical methods. In this paper we develop an analytical model which in principle gives a direct analytical reconstruction formula, but which will be too unstable to use directly. So we will use the analytical model to compute a penalty functional from the noisy measurement data in a parameter space. We do not need optimization algorithms since the space can be fully computed and the minimizer even visually inspected.

This particular problem setting with time-dynamic heat sources \eqref{heatSource} in 3D may not be very common in scientific literature, in 1D similar ideas are used  in the field of material testing using Pulsed Phase Tomography algorithms \cite{BenLardi2010,Marinetti1999}. In our case it is loosely based on the application of breast cancer detection using Fast Fourier Trasform (FFT) as described in \cite{Anbar2001,Button2004} and in \cite{Joro2008a,Joro2008b,Joro2009,Joro2012}. In those papers, by considering the measurement $T(\cdot,t)|_{\partial \Omega}$ a stationary signal and computing its FFT it was hypothized that the frequency spectrum would then reveil abnormal thermodynamic behaviour, for example because of increased blood pulsation and vasodilation. In our paper we build on this idea and investigate how the heat source variation presents itself on the boundary. It turns out that each frequency component $S_n\cos(2\pi f_n+\phi_n)$ of a point heat source results to a temperature signal of the same frequency, but attenuated and phase-shifted according the the frequency.

The main novelty in our approach is the possibility to reconstruct a seemingly complicated heat source with low computational cost, consequently we can plot the complete penalty functional in several dimensions and see how the nature of the inverse problem presents itself there. Another option with potentially less approximation error would be to use Finite Element Method (FEM) or other numerical solvers to compute the forward model for any heat source. Since this can be computationally expensive, to optimize our penalty functional we would use some iterative scheme and hope that the obtained minumum is not local, but global. The present work serves as a proof of concept of the idea of using point source models as building blocks; in the future we aim to investigate linear forward models built from an array of point sources and then the robust inversion, for example using truncated singular value decomposition, to reconstruct source information from boundary temperature measurements.

Unfortunately it would be too much to investigate all of these matters thoroughly in one simple scientific article, and thus we have to make approximations and simplifications, some of which seem severe. We assume that the object is homogenic in material, so that the specific heat capacity, density and the heat conductivity is constant. We don't simulate the IR-sensor technology and their specific noise profiles, different in space and time dimensions, which will always be present in real world measurements, but instead add simulated gaussian (white) noise to the simulated measurements. Also we make use of analytical solutions of the heat equation in the full space, whereas in reality we are always dealing with finite objects and subsequently with certain boundary conditions. Moreover, we do not investigate the actual bio-physical models nor the phenomena they formulate.

We continue with the presentation of our analytical model in section \ref{sec:model}. Then we give a brief review of relevant scientific work in section \ref{sec:review}. After, we describe the reconstruction algorithm in section \ref{sec:algorithm}. The simulation of data is described in section \ref{sec:simulation}. In section \ref{sec:numerics} we present our numerical results, and finally conclude the paper in the final section \ref{sec:conclusions}.

\section{Mathematical model}\label{sec:model}
Define the heat equation by
\begin{equation}\label{HE}
c(x)\frac{\partial T}{\partial t} = \nabla\cdot(\kappa(x)\nabla T) + q(x,t),
\end{equation}
where $c(x)$ is the volumetric heat capacity, $\kappa(x)$ is the heat conductivity, $q(x,t)$ is the volumetric heat source and $T(x,t)$ is the temperature distribution $\R^3\times\R\to\R$. In this document the source has the form \eqref{heatSource} and we first aim to (analytically) compute $T(x,t)$, later in the detection algorithms we aim to compute the location and signal content of $q(x,t)$ from the knowledge of $T(x,t)$. Note that in real life $T(x,t)$ emerges automatically as it is the result of the natural heat diffusion process.


We might not have an analytical solution for the temperature distribution $T(x,t)$ in the most general case. To simplify, we assume that the coefficients $\kappa,c$ and thus the diffusivity $\alpha = \kappa/c$ are constant, and abuse notation a bit to write the heat source
\begin{equation}\label{heatSource2}
q(x,t) = \frac{1}{c}s(t)\delta(x_0),
\end{equation}
and consider
\begin{equation}\label{HE2}
\frac{\partial T}{\partial t} = \alpha\Delta T + q(x,t).
\end{equation}
Then, using standard analytical tools such as the fundamental solution (Gaussian kernel)
\begin{equation}\label{Fundamental}
G(x,t) = \frac{1}{(4\pi\alpha t)^{n/2}}e^{-\frac{x\cdot x}{4\alpha t}},
\end{equation}
where $n=3$ for 3D problems, and Duhamel's principle, we get a solution of \eqref{HE2} with zero initial condition,
\begin{equation}\label{solution3}
T(x,t) = \frac{1}{(4\pi\alpha)^{n/2}} \int_0^t  \frac{e^{-x\cdot x/(4\alpha (t-y))}}{(t-y)^{n/2}}s(y) \dee y,
\end{equation}
where $s(y)$ is the time-varying signal of the heat source \eqref{heatSource2}. We may also change the order of convolution, to get the later useful formula
\begin{eqnarray}
T(x,t) &=& \frac{1}{\pi^{n/2}} \int_0^t  \frac{e^{-x\cdot x/(4\alpha y)}}{(4\alpha y)^{n/2}}s(t-y) \dee y \nonumber \\
&=& \frac{1}{4\alpha}\frac{1}{\pi^{n/2}} \int_0^{4\alpha t}  \frac{e^{-x\cdot x/w}}{w^{n/2}}s(t-\frac{w}{4\alpha}) \dee w \label{solution3b}.
\end{eqnarray}

An important property of heat equations with constant diffusivity $\alpha$, such as \eqref{HE2}, is the following: let $T_1,T_2$ be solutions to
\begin{eqnarray*}
\partial_t T_1 &=& \alpha\Delta T_1 + q_1,\\
\partial_t T_2 &=& \alpha\Delta T_2 + q_2.
\end{eqnarray*}
Then it is clear that $T = T_1+T_2$ solves
$$
\partial_t T = \alpha\Delta T + q_1 + q_2.
$$
This means that by solving analytical solutions $T_n$ for heat sources $q_n$ we also get the solution corresponding to the source $\sum_n q_n$ by a simple summation of the solutions. In particular this principle can be used to combine analytical solutions of multiple separate point sources, and also each heat signal source can be separated into its harmonics, each of which constitutes as a single (time varying) heat source. Moreover this means that given an added time-constant source,
$$
q(x,t) = q_t(x,t) + q_c(x),
$$
the resulting time-dynamic temperature profile can be computed by subtracting the constant temperature solution. 


Of course all of the above mentioned solutions exist only in certain function spaces and with certain theoretical considerations not explained here. What now follows is some extra considerations on how the simple equation \eqref{HE2} differs from more realistic models in section \ref{sec:pennes} and what kind of full-space solution formulas derived from \eqref{solution3} exist in sections \ref{sec:staticSolution} and \ref{sec:dynamicSolution}.

\subsection{About Pennes' bio-heat model}\label{sec:pennes}
According to the well known Pennes' bio-heat transfer model the temperature $T$ in biological tissue satisfies
\begin{equation}\label{pennes}
\rho c \frac{\partial T}{\partial t} = \nabla(k\cdot\nabla T) + w_b\rho_b c_b(T_a-T) + Q_m,
\end{equation}
where $\rho,c$ and $k$ are the density, thermal capacity and heat conductivity of the tissue respectively, $w_b,\rho_b$ and $c_b$ are the perfusion, density and thermal capacity of the blood respectively, $T_a$ is the temperature of arterial blood and $Q_m$ is the metabolic heat source. When compared to \eqref{HE}, the difference is the form of the source term. Of course in any specific situation a model using a source of the general type $q(x,t)$ might be a good approximation, but in general the convection term $w_b\rho_b c_b(T_a-T)$ is needed to model the blood circulation and its effect on the heat dynamics. The Pennes' bio-heat model and its modified versions are well studied, see for example \cite{Ahmadikia2012,Lakhssassi2010,Das2013,Das2014}.

We should also consider finite objects, say domain $\Omega\subset\R^3$ with a boundary $\partial\Omega$, and the full boundary value problem with Neumann-type boundary condition, so that the temperature distribution solves \eqref{pennes} inside $\Omega$ and satisfies
\begin{equation}\label{boundaryCondition}
\kappa\frac{\partial T}{\partial\nu} = h(T_{atm}-T) \quad \textrm{on }\partial\Omega,
\end{equation}
where $h$ is the heat transfer coefficient and $\nu$ is the outwards unit normal vector. This type of boundary condition is usually found in conjunction with the Pennes' model. It is not clear in what way this type of condition will ''affect'' the full-space solution, or in other words how the solution in a finite domain satisfying \eqref{boundaryCondition} differs from the full space solution of \eqref{HE} or \eqref{HE2}.

In the simulations and numerical tests of this paper we use the model \eqref{HE2}; analysis of the error caused by the omission of the convection term and the boundary condition in real world applications is left for future works.

\subsection{Static full-space solution}\label{sec:staticSolution}
Assume a static heat source \eqref{heatSource2} with $s(t) = Q\delta(x_0)$ (and abuse the notation again to infuse the term $c$ into $Q$). We now follow the method described in \cite{Arthur2014,Shi2015}. Take the equation \eqref{HE2} in steady-state and in spherical coordinates, with the origin set at the point source location $x_0$, to get
\begin{equation}\label{HEstatic}
\Delta T = -\frac{Q}{\alpha}\cdot\delta(r),
\end{equation}
where $r$ is the distance from $x_0$. The solution will be spherically symmetric. The Laplace-operator simplifies to only include $r$-terms, and we get
\begin{equation*}
\frac{1}{r^2}\frac{\partial}{\partial r}\left(r^2\frac{\partial T}{\partial r}\right) = -\frac{Q}{\alpha}\cdot\delta(r).
\end{equation*}
This differential equation is easily solved for its general solution 
\begin{equation}\label{staticSolution1}
T(r) = \frac{Q}{4\pi\alpha r} + C,
\end{equation}
where $C$ is a constant. In \cite{Arthur2014} this analysis is followed by a method to then compute the unknown $Q$ and the depth $d$, defined as the shortest distance between $x_0$ and the boundary, from the boundary temperature distribution $T(x)|_{\partial \Omega}$. In \cite{Shi2015} another method is presented, using curve fitting to plot the ''Q-r'' -curve which according to their clinical study has different profile for different diagnostic cases. Neither of these methods are using the usual, regularized mathematical machinery from the field of inverse problems, and so we continue with our own implementation in section \ref{sec:algorithm}. One aspect we take care of is to not only search for the depth of the heat source, but the location as well, since it might not be clear where the point of maximum temperature is located exactly in a noisy measurement.

It should be noted here that in \cite{Arthur2014} there is also a step to filter the temperature distribution to account for non-singular heat sources. We omit this step, and just simply test our resulting algorithm in a simulation where the heat source is not as small as possible, but a decently sized sphere.

\begin{figure}
\begin{picture}(120,90)
\epsfxsize=10cm
\put(10,0) {\epsffile{\imagepath 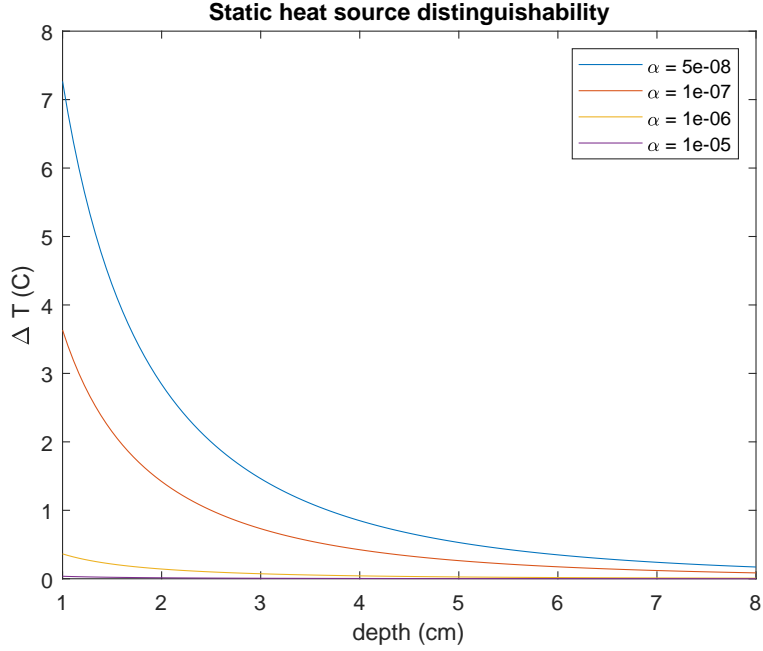}}
\end{picture}
\caption{\label{fig:staticDT} Distinguishability of a static, spherical heat source with $Q = 29000 W/m^3$ and radius $R=1 cm$. Diffusivity $\alpha$ is in units $m^2/s$. See text for further explanation on distinguishability.}
\end{figure}

In figure \ref{fig:staticDT} we plot the distinguishability of an abnormal heat source using different values of diffusivity. The heat source is spherical with values $Q = 29000 W/m^3$ and radius $R=1 cm$, mimiking the metabolic rate of a cancerous tumor. The distinguishability here means the difference of temperature between the hottest spot and 5cm away on a plane. In real applications the distinguishability will be less because of noise and boundary conditions such as \eqref{boundaryCondition}. Moreover the distinguishability here is a crude approximation of what it would be like to actually perform an analysis on the boundary temperature data, which entails the full temperature data between the hottest spot and the nearby area, up to 5cm radius. Fat and muscle tissue have diffusivity of 0.98E-7 and 1.31E-7 respectively, which will be suitable for the static detection of heat sources, based on figure \ref{fig:staticDT}. In the next section we will consider dynamic heat sources which have larger optimal diffusivities depending on the frequency of the source.

\subsection{Dynamic full-space solution}\label{sec:dynamicSolution}
The idea of what follows can be found in the book \cite{Carslaw1959} and is also available in the lecture note \cite{M445}. Assume the time-variable heat source \eqref{heatSource2} with one harmonic component,
\begin{equation*}
q(x,t) = S_n\cos(2\pi f_n t + \phi_n)\delta(x_0).
\end{equation*}
Then, by writing 
$$
f(w) = \frac{e^{-x\cdot x / w}}{w^{n/2}}, \quad \lambda = 2\pi f_n t + \phi_n,\quad \beta = \frac{2\pi f_n}{4\alpha} \quad A = \frac{S_n}{4\alpha\pi^{n/2}},
$$
from \eqref{solution3b} and trigonometric identities, as $t\to\infty$, we obtain
\begin{eqnarray*}
T(x,t) &\approx& A\left(\cos(\lambda)\int_0^\infty\cos(\beta w)f(w)\dee w + \sin(\lambda)\int_0^\infty\sin(\beta w)f(w)\dee w \right) \\
&=& A\left(\cos(\lambda)\Re(\hat{f}(\beta)) + \sin(\lambda)\Im(\hat{f}(\beta))\right),
\end{eqnarray*}
where $\hat{f}$ is the half-line Fourier-transform of $f$. In 3D we have
$$
\hat{f}(\omega) = (\sqrt{2+\sqrt{2}}+\rmi \sqrt{2-\sqrt{2}})\frac{\sqrt[4]{\omega}}{\sqrt{\abs{x}}}K_{1/2}(2\sqrt{-\rmi\omega}\abs{x}),
$$
where $K_{1/2}$ is a modified Bessel function
$$
K_{1/2}(z) = \sqrt{\frac{\pi}{2}}\frac{e^{-z}}{\sqrt{z}}.
$$
Thus we get
$$
\hat{f}(\omega) = \frac{\sqrt{\pi}}{\sqrt{2}\abs{x}}(1+\rmi)e^{-\sqrt{2\omega}\abs{x}}\left(\cos(\sqrt{2\omega}\abs{x})+\rmi\sin(\sqrt{2\omega}\abs{x})\right),
$$
and subsequently, using again some trigonometric identities,
\begin{eqnarray}
T(x,t) &=&  \frac{A\sqrt{\pi}}{\sqrt{2}\abs{x}}e^{-\sqrt{2\beta}\abs{x}}\left(\cos(\lambda - \sqrt{2\beta}\abs{x}) + \sin(\lambda - \sqrt{2\beta}\abs{x})\right) \nonumber \\
 &=& \frac{A\sqrt{\pi}}{\abs{x}}e^{-\sqrt{2\beta}\abs{x}}\cos(\lambda - \sqrt{2\beta}\abs{x} - \pi/4) \nonumber \\
 &=& \frac{S_n}{4\alpha\pi\abs{x}}e^{-\sqrt{\frac{\pi f_n}{\alpha}}\abs{x}}\cos(2\pi f_n t + \phi_n - \sqrt{\frac{\pi f_n}{\alpha}}\abs{x} - \pi/4) \label{timeSol}. 
\end{eqnarray}
Note that the original heat signal source will inhibit phase change and amplitude attenuation depending on the frequency and distance from the source; in particular the amplitude attenuation is exponential.

Take the full heat source of  \eqref{heatSource}. By using the summation property and \eqref{timeSol}, the corresponding analytical solution of the temperature distribution is
\begin{equation}\label{multiSol}
T(x,t) \approx \frac{1}{4\alpha\pi\abs{x}}\sum_{n=1}^N \tilde{S}_n\cos(2\pi f_n t +\tilde{\phi_n}),
\end{equation}
where
\begin{eqnarray}
\tilde{S}_n &=& S_n e^{-\sqrt{\frac{\pi f_n}{\alpha}}\abs{x}}, \\
\tilde{\phi}_n &=& \phi_n - \sqrt{\frac{\pi f_n}{\alpha}}\abs{x}-\pi/4.
\end{eqnarray}
Many specific heat sources can be divided into its components in a similar way, and their respective temperature distributions deduced using the harmonics. Some of these are presented in the next section.

\begin{figure}
\begin{picture}(120,90)
\epsfxsize=10cm
\put(10,0) {\epsffile{\imagepath 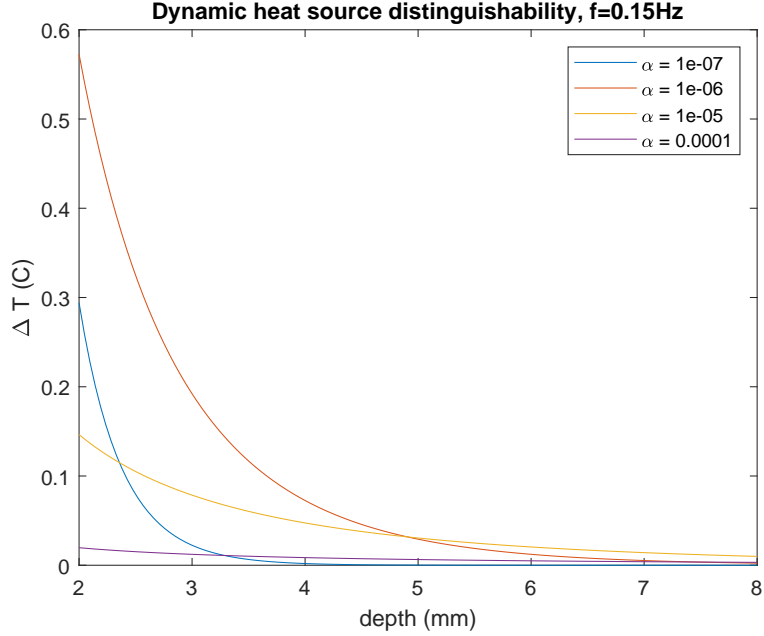}}
\end{picture}
\caption{\label{fig:dynamicDT1} Distinguishability of a dynamic, spherical heat source with $S = 29000 W/m^3$, radius 1 cm and frequency 0.15 Hz. Diffusivity $\alpha$ is in units $m^2/s$. See text for further explanation on distinguishability.}
\end{figure}

\begin{figure}
\begin{picture}(120,90)
\epsfxsize=10cm
\put(10,0) {\epsffile{\imagepath 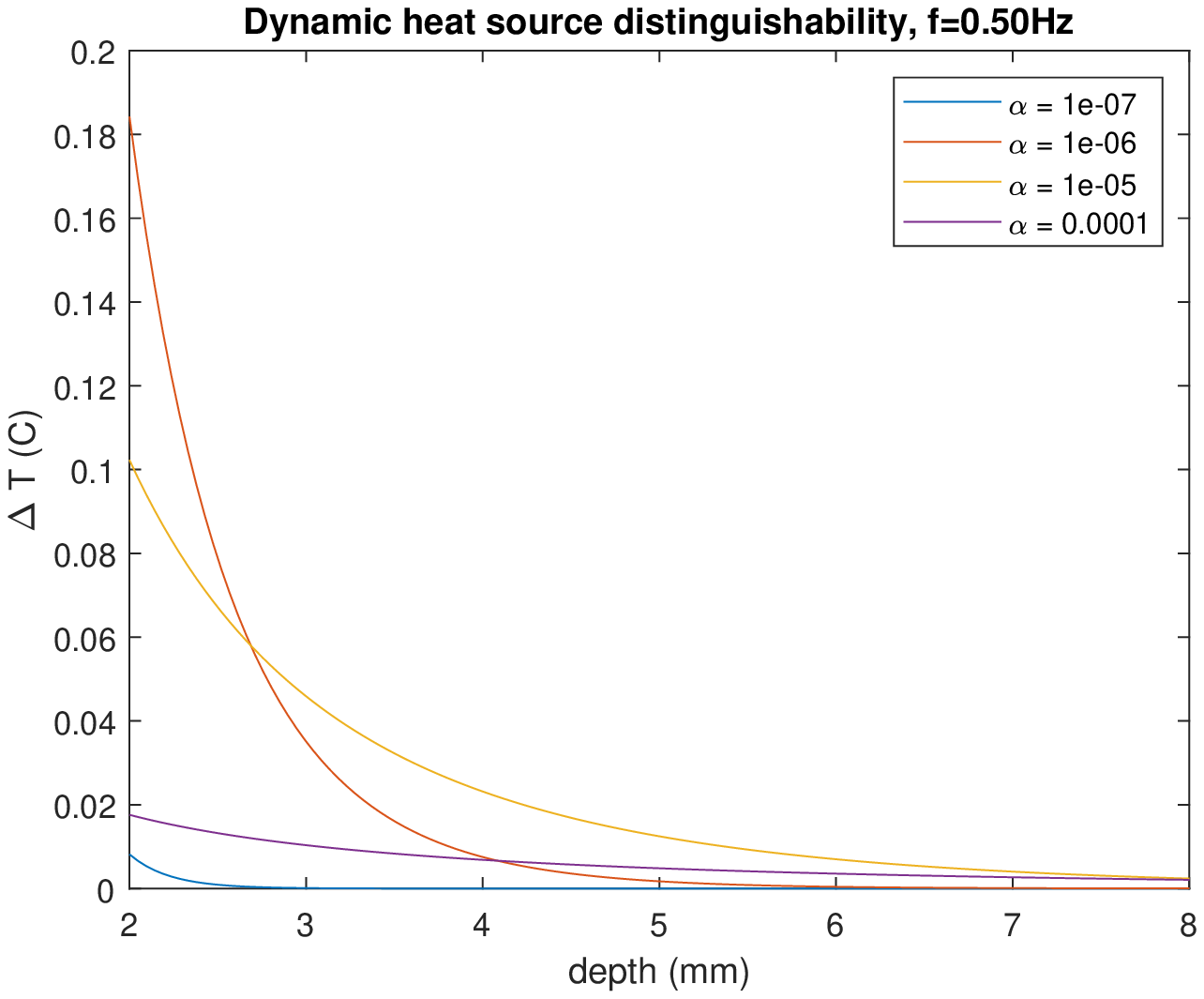}}
\end{picture}
\caption{\label{fig:dynamicDT2} Distinguishability of a dynamic, spherical heat source with $S = 29000 W/m^3$, radius 1 cm and frequency 0.50 Hz. Diffusivity $\alpha$ is in units $m^2/s$. See text for further explanation on distinguishability.}
\end{figure}

\begin{figure}
\begin{picture}(120,90)
\epsfxsize=10cm
\put(10,0) {\epsffile{\imagepath 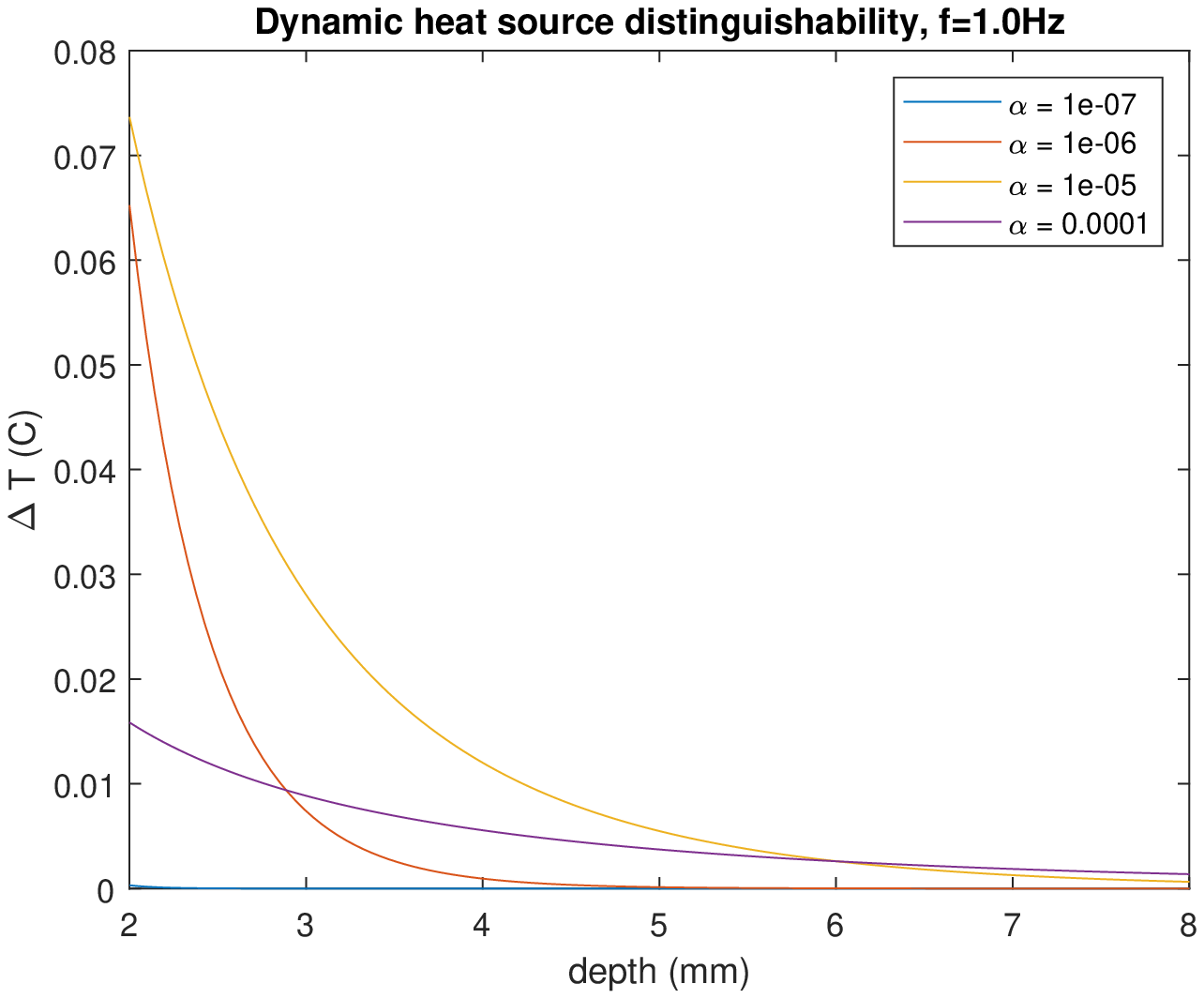}}
\end{picture}
\caption{\label{fig:dynamicDT3} Distinguishability of a dynamic, spherical heat source with $S = 29000 W/m^3$, radius 1 cm and frequency 1.0 Hz. Diffusivity $\alpha$ is in units $m^2/s$. See text for further explanation on distinguishability.}
\end{figure}

Note the different attenuation terms in \eqref{timeSol}. The term $1/\abs{x}$ is the same for all frequencies as in the static case \eqref{staticSolution1}. The exponential term depends also on the frequency. Both terms depend on the diffusivity. As in figure \ref{fig:staticDT}, it is interesting to see the distinguishability for different diffusivity coefficients and frequencies. The scale of the depth and diffusivity in figures \ref{fig:dynamicDT1}, \ref{fig:dynamicDT2} and \ref{fig:dynamicDT3} are different from the static case because the attenuation is much more severe because of the exponential term. In this dynamic case, it is more difficult to estimate what the distinguishability, computed in the same way as in the static case, means. In principle the amplitude information is computed using FFT and if the phenomena are perodic and we have long enough acquisition time, we are able to detect very small temperature variation amplitudes.

\subsection{Some examples}
Consider a phase-modulated source and its series expansion
\begin{equation}\label{modulation}
s(t) = \cos(\omega t + B\sin(\omega_m t)) =  \sum_{k=-\infty}^{\infty}J_k(B)\cos((\omega+k\omega_m)t),
\end{equation}
where $\omega_m$ is the frequency of modulation and $J_k(x)$ is a Bessel function of the first kind, of order $k$. Following earlier reasoning, each component in the sum \eqref{modulation} contributes to a solution given by \eqref{timeSol}; each component will be attenuated  and phase-shifted according to their frequency, higher frequencies more than lower. The source term itself has main frequency $\omega$ and symmetrical sidelobes of frequencies $n\omega_m$ with amplitudes $J_n(B)$. The smaller the modulating amplitude $B$ is, the smaller these sidelobes are; when $B\geq 1$ the sidelobes become significant and will tail off as a sinusoid in their amplitude. The important thing is that for larger amplitude, say $B=2$, the amplitude of the main frequency will be lower than the first sidelobes, and in general for even larger $B$ the largest harmonic components are (symmetrically) further away from $\omega_c$. Thus without knowing the attenuation behaviour described by equation \eqref{timeSol}, we could register a lower frequency from the temperature variation than the actual source frequency.



Consider an amplitude modulated source and its trigonometric expansion
\begin{eqnarray*}
s(t) &=&  (1+M\cos(\omega_a t))\cos(\omega t + \phi) \\
&=& \cos(\omega t + \phi) + \frac{M}{2}[\cos((\omega + \omega_a)t + \phi) + \cos((\omega - \omega_a)t)].
\end{eqnarray*}
The first term can be thought of as our signal of interest, the second introduces sidelobes situated $\omega_a$ away from the original, with amplitudes $M/2$. Again as with the phase-modulated source, each of these components contributes to a solution given by \eqref{timeSol}.


\section{Literature review}\label{sec:review}
The heat source reconstruction problem has been studied in many aspects: in spatial dimensions 1,2 and 3, with various forms of the heat source, with constant and variable heat conductivity and with various assumptions on the measurement, both Dirichlet and Neumann conditions on either partial or full boundary. Many papers exist on the stability of the problem, that is an estimate on how the error between the reconstructed solution and the true solution behaves as a function of the measurement error. Some papers concentrate on the mathematical theory, others include some numerical tests and even numerical algorithms.

In \cite{Dalir2014} analytical solutions of the three dimensional heat equation, using a series expansion, are computed for a multi-layered sphere. In \cite{Flint2018} analytical solutions are used to model moving heat sources i.e. welding process. Both of these papers could be the inspiration for an even more advanced analytic modeling and subsequent heat source detection.

The papers \cite{Davoodi2011,Liao2009,Liu2010} deal with one dimensional numerical heat source reconstruction, but they do not have any regularization even though the problem is ill-posed. A more mathematical background  dealing with uniqueness and establishing the problem as an inverse problem is given by J.~R.~Cannon in \cite{Cannon1968,Cannon1982} and in numerous other papers by the same author. Continuing the analytical approach we have the theoretical stability analysis of papers \cite{Cannon1976,Cannon1986,Cannon1991,Talenti1982,Trong2010,Yamamoto1993}. Examples of regularization strategies to battle against the ill-posedness are found in \cite{Cheng2010,Cheng2012,Dou2009,Dou2012,Hao2017,Wang2004}. 

Finally, in \cite{Renault2007} we have a regularized, numerical algorithm to reconstruct time-dynamic heat source in 2D, for the purposes of material testing and the analysis of thermomechanical effects. The method is further developed in \cite{Renault2010}. 

Please see also the references in the publications above, as the list presented here is only a very limited and somewhat arbitrary subset of all of the work done in this field.

\section{Algorithms}\label{sec:algorithm}
We can think of several classes of heat source reconstruction algorithms, from the most general case to special cases using a priori information. In order to describe the algorithms we need more notation. The boundary of the object is $\partial\Omega = B(x)$. The ''full'' measurement is $M(x,t) = T(x,t)|_{B(x)}$, where $t\in[0,t_{max}]$. From $M$ we may manually choose a more interesting neighbourhood $M_0(x,t)$ where $x\in \mathcal{N}\subset B$, containing a visible hot spot we want to analyze in more detail. The FFT of these measurements are $\hat{M}$ and $\hat{M}_0$, for each spatial point $x$, more specifically we write
\begin{eqnarray}
\hat{M} &=& \{A_n,\theta_n\},\quad n=1,\ldots,N,\\
\hat{M}_0 &=& \{A^0_n,\theta^0_n\},\quad n=1,\ldots,N,
\end{eqnarray}
where the amplitudes $A_n,A^0_n$ and phases $\theta_n,\theta^0_n$ for each $n$ describe the signal component corresponding to the frequency $f_n = n/N_s\cdot f_s$, where $f_s$ is the sampling rate and $N_s$ is the number of samples. 
The index $n=1,\ldots,N$ is meant to go through the interesting frequencies, if not all available from FFT, in which case we would have roughly $N=N_s/2$ because of Shannon-Nyquist sampling theorem.

For each point $x_0$ and frequency $f_n$ we can write
\begin{eqnarray}
F(x) &=& \frac{1}{4\pi\alpha\abs{x_0-B(x)}},\\
D_n(x) &=& \sqrt{\frac{\pi f_n}{\alpha}}\abs{x_0-B(x)},
\end{eqnarray}
and subsequently we can express the diffused source component of \eqref{timeSol} as
\begin{equation}\label{modelSolution2}
T(x,t) = \tilde{S}_n(x)\cos(2\pi f_n t + \tilde{\phi}_n(x)),
\end{equation}
where
\begin{eqnarray}
\tilde{S}_n(x) &=& F(x)e^{-D_n(x)}S_n,\\
\tilde{\phi}_n(x) &=& \phi_n - D_n(x) - \pi/4.
\end{eqnarray}
Also the static solution \eqref{staticSolution1} can then be written as
\begin{equation}\label{modelSolution1}
T(x) = QF(x) + C.
\end{equation}

There are many ways to arrange the detection algorithm in terms of whether it is model or data driven. In our algorithms, we take a grid of source location candidates $x_0$ and output the one having the smallest penalty functional. But to compute the penalty, we first use the model and the data together to compute some intermediate steps. In the static case, we compute a maximum and minimum temperature value on $M_0$ from which we get the model complient parameters $Q,C$ for that particular source location. Then, the penalty functional is
\begin{equation}\label{staticError}
\epsilon = \norm{M_0(x)-(QF(x)+C)}_{L^2(B)}.
\end{equation}

In the dynamic case, for each candidate point source $x_0$ we compute its nearest boundary point $C$ on $M_0$ and then the values $\{A^0_k,\theta^0_k\}$ for each frequency component $k$. Using the dynamic model \eqref{modelSolution1} these are then turned to reconstructions $\{ S_k,\phi_k\}$. However, we both cut away impossible spectral content and find the best fit from the source locations. By defining truncation parameters $A_t$ and $M_t$, we define the set $M_0^k\subset M_0$ as the set for which has spectral amplitude over $A_t$ (cutting off too small data) and which has the size of at least $M_t$ (cutting off too little pieces of measurement data). The lower $A_t$ and higher $M_t$ is, the more robust we need elsewhere in the algorithm. For each source candidate, for the amplitude data, we use the penalty functional 
\begin{equation}\label{dynamicError1}
\epsilon_k^A = \norm{A_k-\tilde{S}_k }_{L^2(M)},
\end{equation}
where the profile $\tilde{S}_k$ is given by the source candidate $x_0$, the computed value $A^0_k$ and the dynamic model \eqref{modelSolution1}. The total amplitude penalty functional is
\begin{equation}\label{dynamicError2}
\Sigma_1 = \sum_{k=1}^N\epsilon_k^A,\quad M_0^k\neq\emptyset.
\end{equation}

For the phase penalty functional, we use
\begin{equation}\label{dynamicError3}
\epsilon_k^\theta = \norm{\theta_k-\tilde{\phi}_k^1 }_{L^2(M_0^k)},
\end{equation}
where $\tilde{\phi}_k^1$ is the model computed profile (at each source location $x_0$) normalized to have the same scale as the measurement $\theta_k$ on $M_0^k$. Similarly to the amplitude penalty we write
\begin{equation}\label{dynamicError4}
\Sigma_2 = \sum_{k=1}^N\epsilon_k^\theta,\quad M_0^k\neq\emptyset.
\end{equation}
We were not able to combine the information given by the amplitude and phase profiles in a good way, so instead we get two reconstructions using either information.

Lastly, we denote the normalized collection of harmonics by $\{S^1_k,\phi_k\}$, satisfying $S_k = A S_k^1$ for each $k=1,\ldots,N$, where $A$ is the amplitude of the signal with harmonic content $\{S^1_k,\phi_k\}$.


\subsection{General algorithm descriptions}
There are many possibilities how to arrange the algorithm to find the location $x_0$ and the spectral content $\{S_n,\phi_n\}$ from the measurement, depending also on what we expect the situation to be. We can think of the following examples:
\begin{itemize}
\item[(A1)] Given $M$, reconstruct $q(x,t)$.
\item[(A2)] Given $M$, reconstruct a collection of point sources $\{x_0\}_m$ with signals $\{S_n,\phi_n\}_m$.
\item[(A3a)] Given $M_0$, reconstruct $Q,C,x_0$.
\item[(A3b)] Given $M_0$, reconstruct $x_0$ and $\{S_n,\phi_n\}$.
\item[(A4)] Given $M_0$ and $\{S^1_n,\phi_n\}$, reconstruct $x_0$ and $A$.
\end{itemize}
The algorithm (A1) will be difficult to realize and is not in the scope of this paper. Our model only considers point-like sources. In principle by using the linearity of the heat equation with constant diffusion constant one might be able to approximate different geometries of heat sources by using a collection of point sources. A moving pulse-like heat source can be approximated by a series of point sources, but then an issue of phase cancellation emerges, when two different points sources produce two temperature variations in different phase resulting to reduced temperature signal. Part of the algorithm (A1) would be to reconstruct the geometry of moving sources from the phase information.

To understand the difficulty of algorithm (A2), consider the fact that each point source will be visible in the part of $M$ closest to them, as a hot spot from which we want reconstruct the parameters of the point source. The algorithm (A2) would need an automatic scan of the full measurement $M$ in such a way that these visible hot spots are first scanned, after which they are analyzed separately. Again, we do not attempt this, but instead assume that the user is able to choose an interesting neighbourhood which contains a hot spot, to reconstruct each heat source $x_0$ of interest.

The algorithm (A3a) uses the static model \eqref{modelSolution1} and is described in more detail in the next section. The algorithm (A3b) uses the dynamic model \eqref{modelSolution2} and is described in section \ref{dynamicAlgorithm}.

The algorithm (A4) is thought of as a simple way of using a priori information about the signal, for example we could have another means of measuring the heat source signal content in conjunction with the temperature profile on the boundary after which the algorithm would search for this signal content in the temperature data. The details for this is left in future works.

\subsection{Source power reconstruction}
Note that a finite domain has a boundary and boundary conditions. If there are positive heat sources and perfect insulation at the boundary, so that the domain would not lose heat to the surroundings, then the temperature of the domain would increase infinitely. A more realistic scenario, which would be true for all objects in the physical world, is a boundary condition such as equation \eqref{boundaryCondition} which would result to a stable state after a sufficient stabilization time. The stable temperature of the domain would be dependent on the values of the heat source power and the heat transfer coefficient on the boundary. Having a larger power heat source and larger heat transfer coefficient could result to approximately the same average temperature as smaller power heat source and smaller heat transfer coefficient; thus our analytical model which takes none of this into account can never perform an absolute reconstruction of the heat power.

However, using the equation \eqref{modelSolution1} we may still compute an absolute value for the power amplitude $Q$, even if it won't be a physically correct value. Subsequently our model can still compare different powers of heat sources since they can be reconstructed having different values of $Q$. To do this we use the following sub-algorithm (AQ) for each candidate location $x_0$:
\begin{enumerate}
\item Sort the values of $M$ or $M_0$ and take an average of the first and last $\mathrm{Nav}$ values to get the minimum temperature $T_1$ and maximum temperature $T_2$ of the measurement patch. Besides the sorted vector of measurements, use their locations to compute the corresponding average locations $x_1$ and $x_2$ and their distances $r_1$,$r_2$ from $x_0$.
\item Since we seek for values $Q,C$ for the model \eqref{modelSolution1}, from there we get the formula
$$
Q = (T1-T2)(1/(4\pi\alpha r_1) - 1/(4\pi\alpha r_2))^{-1}.
$$
\item Compute the values $C_1,C_2$ again from the equation \eqref{modelSolution1} using $T_1,T_2,r_1,r_2$ and $Q$, and finally compute
$$
C = (C_1+C_2)/2.
$$
\end{enumerate}
In this way we can compute for any candidate $x_0$ the power $Q$ and the parameter $C$. If $M$ or $M_0$ is very noisy the above will perform less well, but the integer $\mathrm{Nav}$ can be used to alleviate this problem a little.

\subsection{Source spectrum reconstruction}
The case of a dynamic heat source becomes more complicated as we have one more dimension in data and source parameters. The biggest problem is the sensible use of the phase information; the model could easily predict a total phase change of $6\pi$, for example, but the FFT processing of any simulated or real data will not be accurate throughout the measurement patch, but instead only inside a small radius. In addition to this the overall level of the phase error norm will be different form the amplitude level, thus making it difficult to combine the phase and amplitude information. The algorithm presented here is only one somewhat simple suggestion and should be investigated more in a subsequent paper.

One of the features is that we automatically cut away frequency information that is not large enough and only reconstruct the best source fit for the remaining data. In some future work when tackling more noisy data this approach will probably be too limiting. Other feature that we had to reluctantly settle with is that the algorithm finds the best fit separately using amplitude and phase information since we did not find a good way to combine the differently behaving penalty functionals.

In the dynamic case the measurement $M$ or $M_0$ has been processed by FFT to get the collection $\hat{M}$ or $\hat{M_0}$ of amplitudes and phases. Define an amplitude limit $A_t$ below which we discard amplitude data. Then, given any source location candidate $x_0$ we use the following sub-algorithm (AA):
\begin{enumerate}
\item Compute the boundary point $C$ closest to the source and write $A_n^C$ for the amplitude values on a neighbourhood of $C$, for each frequency bin $n$. Write $r=\abs{x-x_0}$.
\item For each $n$, use \eqref{modelSolution2} and $A_n^C$ to compute the model predicted amplitude profile $\tilde{A}_n = A_n^CF(r)e^{-D_n(r)}$. Compute the sub-measurement patch $\tilde{M}_n = \{ x\in B(x):\quad \tilde{A}_n>A_t\}$. 
\item If $\tilde{M}_n = \emptyset$ or a very small set, put $S_n = \phi_n = 0$ as the reconstruction. Otherwise use $\tilde{M}_n$ as the domain where the phase error is computed.
\end{enumerate}

The phase profile $\phi_n$ needs special processing as its raw values are always inside $[0,2\pi)$. We use Matlab's unwrap-function to first make the phase profile continuous, then put the phase value on $C$ to zero by subtracting the value computed earlier in conjunction with the amplitudes $A_n^C$. The analytical model will produce a reference phase profile that is already continuous, and we process it the same way so that on the boundary point $C$ the phase will be zero. The phase information can be averaged by first changing to unit vectors in 2D, then averaging the vectors and finally changing back to the angle of the averaged vector.

Since we only use $\tilde{M}_n$ to compute the phase error, two problems arise. First, the smaller the area of $\tilde{M}_n$, the smaller the magnitude of the phase error. Second, the phase information outside of $\tilde{M}_n$ will not contribute to the error, which is clearly not optimal. The first issue is somewhat solved by dividing the error norm by the area of $\tilde{M}_n$. The second issue is not properly addressed, but instead we compute the optimal source location and the accompanying spectral content using the amplitude and the phase information separately.

\subsection{Static reconstruction algorithm}
The algorithm (3a) in more detail is the following:
\begin{enumerate}
\item From $M$, pick an interesting neighbourhood $M_0$ manually.
\item Take a collection of source candidates $x_0$ and for each of them use (AQ) to compute the parameters $Q$ and $C$. For each triplet $(x_0,Q,C)$ compute $\epsilon$ from \eqref{staticError}.
\item The reconstruction is the triplet having the smallest error $\epsilon$.
\end{enumerate}

\subsection{Dynamic reconstruction algorithm}\label{dynamicAlgorithm}
The algorithm (3b) in more detail is the following:
\begin{enumerate}
\item From $M$, pick an interesting neighbourhood $M_0$ manually.
\item Compute FFT of $M_0$ to get the amplitude and phase data.
\item Take a grid of source candidates $x_0$ and for each use the sub-algorithm (AAn) to get the model-conforming amplitudes $A_n^C$ and usable phase areas $\tilde{M}_n$.
\item For each $x_0$ and each $k$ that is not cut away, compute amplitude and phase penalty functionals\eqref{dynamicError1} and \eqref{dynamicError3}.
\item The reconstruction is the location $x_0$ and the accompanying signal $\{S_k,\phi_k\}$ using either \eqref{dynamicError2} or \eqref{dynamicError4}.
\end{enumerate}

\section{Simulation}\label{sec:simulation}
For medical imaging, ideally we would want to simulate data on a realistic 3D geometry, using realistic values of heat conductivity and capacity for anatomically correct tissues. The discrete computing mesh should enable high accuracy in the parts of the domain where we measure physical properties, but should be adaptively course to keep down computational demands. Unfortunately in this paper we needed to keep things simple. The domain is a half ball of radius 100, the computational mesh is similar everywhere in the domain and has maximum edge length of 10 for the static examples, 14 for dynamic examples, the diffusivity is constant 2 and the values of heat generation and heat conduction outwards are chosen so that we get somewhat realistic (between 20 and 40 degree Celsius) temperature values inside the object.

We investigate static and dynamic ball-shaped heat sources of radius 10. For the static cases, there are three different locations of the heat source each with four different volumetric power. For the dynamic case we only compute two locations closest to the boundary and two even larger amplitudes of power modulation compared to the static case, because the effect of attenuation is severe.

We use Finite Element Method (FEM) to compute the temperature profile in three spatial dimensions in addition to the time dimension. For the static cases we only use the temperature profile at the last available discrete time position. For the dynamic cases we use a time-dynamic profile of the last ten seconds of the simulation. In general we computed several rounds of FEM solutions in a serial manner, so that the temperature profile at the end of one simulation was used as the initial value for the next simulation. In this way we could use a course time dimensional computation for the stabilization period, and finer time dimensional computation at the very end  after the system had stabilized.

The dynamic solution is transformed to amplitude and phase information using FFT with frame rate of 10Hz without windowing since it was not needed in our particular choices of frequency and acquisition time. Before FFT we also removed linear trend of warming from the simulated measurement data.

Accuracy problems were encountered because of the four dimensional computations and subsequent need for course FEM mesh in spatial dimensions. We interpolated the temperature profile at radius 90 in the static examples and at radius 85 in the dynamic examples, inside the domain, as the boundary measurement. The dynamic simulations were problematic especially since the dynamic effects are exponentially attenuated, and thus very small in general. 

We simulate noise by adding Gaussian white noise at several levels. One percent noise means the variance of the Gaussian distribution is one percent of the maximum temperature difference in the measurement patch $M_0$.

\begin{figure}
\begin{picture}(120,110)
\epsfxsize=10cm
\put(10,0) {\epsffile{\imagepath 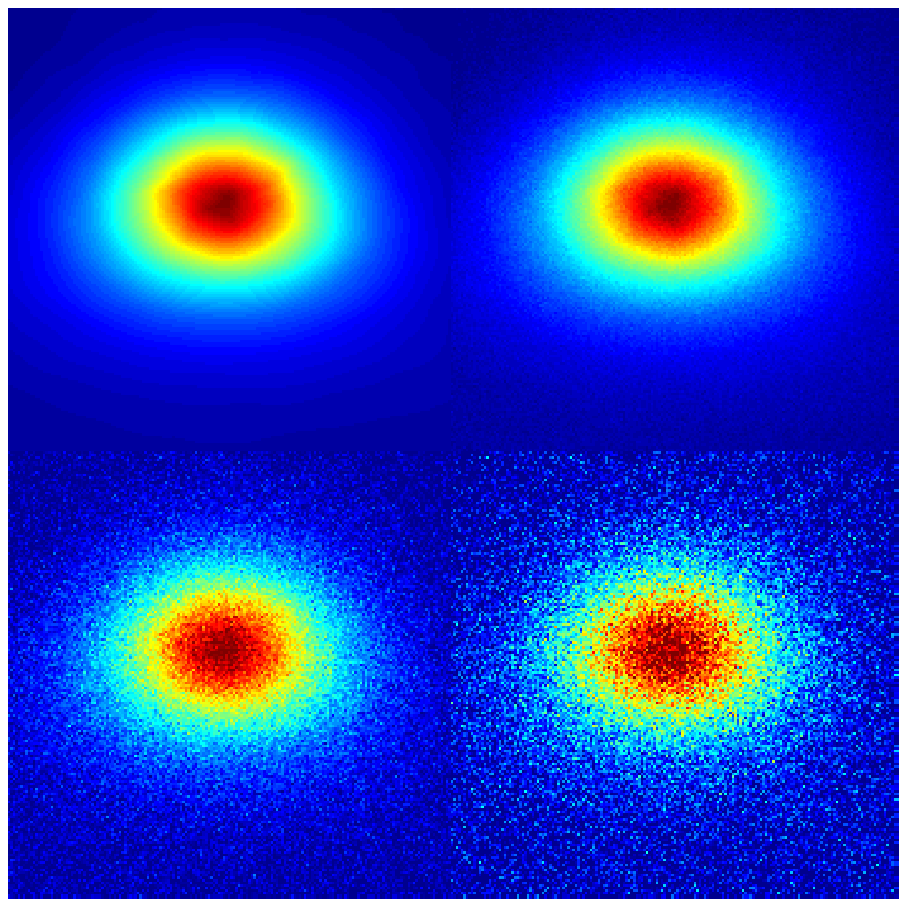}}
\end{picture}
\caption{\label{fig:noiseExample1} The temperature measurement patch of heat source A with power 1.0, non-noisy simulation and three simulated noise levels.}
\end{figure}

\begin{figure}
\begin{picture}(120,110)
\epsfxsize=10cm
\put(10,0) {\epsffile{\imagepath 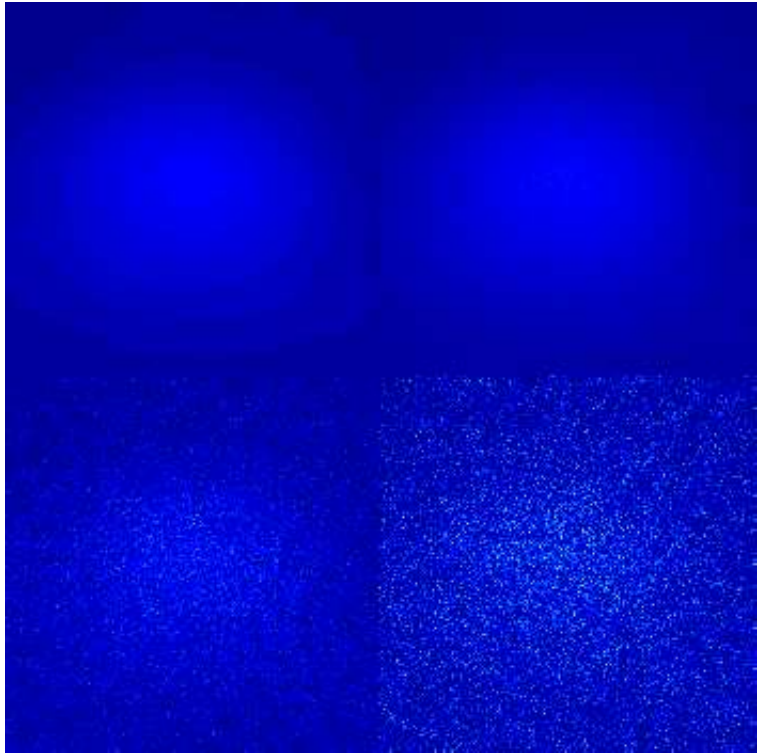}}
\end{picture}
\caption{\label{fig:noiseExample2} The temperature measurement patch of heat source C with power 0.8, non-noisy simulation and three simulated noise levels.}
\end{figure}

The parameters of the heat sources are given in table \ref{parameterTable}. In figures \ref{fig:noiseExample1} and \ref{fig:noiseExample2} we have an example of manually picked temperature measurement patch $M_0$ near the boundary, non-noisy and with three different simulated noise levels, from static cases A1 and C2. The two figures use the same temperature scale and colormap and shows how much smaller and more difficult to detect the ''hot spot'' becomes as the heat source is located deeper in the domain and has a smaller power.

\begin{table}
\begin{tabular}{c|c|c}

             & static, R=10                        & dynamic, f=0.2Hz \\
center       & $Q_1 \quad Q_2 \quad Q_3 \quad Q_4$ & $S_1 \quad S_2$ \\
\hline
A [40 40 50] & $1.0 \quad 0.8 \quad 0.6 \quad 0.4$ & $2.2 \quad 1.8$ \\
B [35 35 45] & $1.0 \quad 0.8 \quad 0.6 \quad 0.4$ & $2.2 \quad 1.8$\\
C [30 30 40] & $1.0 \quad 0.8 \quad 0.6 \quad 0.4$ &                                     \\

\end{tabular}
\caption{\label{parameterTable} Heat source parameters for the simulations. }
\end{table}

\section{Numerical tests}\label{sec:numerics}
We used truncation values $A_t=0.02$ and $M_t = 900$ pixels. In the static case we took $40^3$ different source candidate locations, 40 in each spherical dimension to cover the whole domain $\Omega$ near the measurement patch $M$. In the dynamic case we used $30^3$ different source candidate locations, giving inherently a bit smaller detection resolution.
\subsection{Static reconstructions}
\begin{figure}
\begin{picture}(120,160)
\epsfxsize=10cm
\put(10,0) {\epsffile{\imagepath 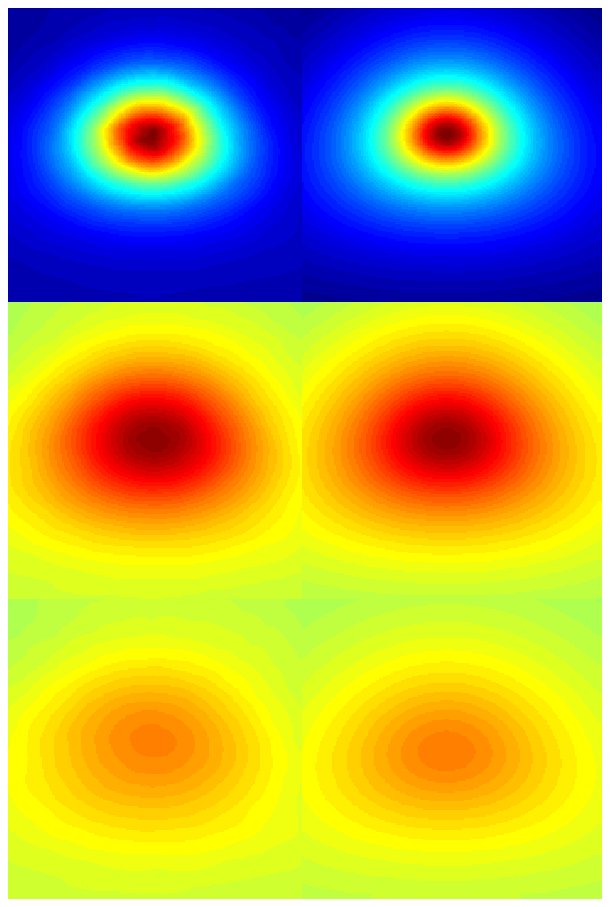}}
\end{picture}
\caption{\label{fig:reconExample} On the left: the measurement patches for heat sources A,B and C with power 1.0. On the right: temperature patches that gave the lowest error.}
\end{figure}

In table \ref{results1} are the results of the static simulations. The location column has the distance from the reconstructed source center to the true source center in same units as the simulations. The power columns give the reconstructed power compared to the non-noisy reconstruction of source A1, which should be the easiest to reconstruct. In figure \ref{fig:reconExample} we see example reconstructions, or the best fit temperature profiles, to show how close the point source model computed measurements are to the simulated measurements.

\begin{table}
\begin{tabular}{c|c|cccc}
            &  location error & $Q_1$ & $Q_2$ & $Q_3$ &  $Q_4$  \\
\hline
center A    &       &                 &       &       &        \\
no noise    &  5.0  &            1    &  0.83 &  0.56 &  0.39  \\
noise 1\%   &  5.2  &          1.1    &   1.0 &  0.79 &  0.48  \\
noise 5\%   &  10.2 &          1.6    &   1.5 &   1.1 &  0.67  \\
noise 10\%  &  10.5 &          1.4    &   1.3 &   1.0 &  0.76  \\
\hline
center B    &       &                 &       &       &        \\
no noise    & 6.6   &           0.90  &  0.48 &  0.38 &  0.34  \\
noise 1\%   & 7.0   &           1.0   &  0.68 &  0.52 &  0.57  \\
noise 5\%   & 19    &           1.0   &  0.68 &  0.87 &  0.52  \\
noise 10\%  & 14    &           2.7   &  0.80 &  0.66 &  0.75  \\
\hline
center C    &       &                 &       &       &        \\
no noise    & 11    &           0.61  &  0.41 &  0.33 &  0.34  \\
noise 1\%   & 8.7   &           0.90  &  0.56 &  0.52 &  0.80  \\
noise 5\%   & 19    &           1.6   &  0.50 &  0.59 &  0.44  \\
noise 10\%  & 19    &           3.0   &  0.84 &  0.68 &  0.69  
\end{tabular}
\caption{\label{results1} Static reconstructions for the three source locations, four power levels and three simulated noise levels.}
\end{table}

\subsection{Dynamic reconstructions}
\begin{figure}
\begin{picture}(120,100)
\epsfxsize=12cm
\put(0,0) {\epsffile{\imagepath 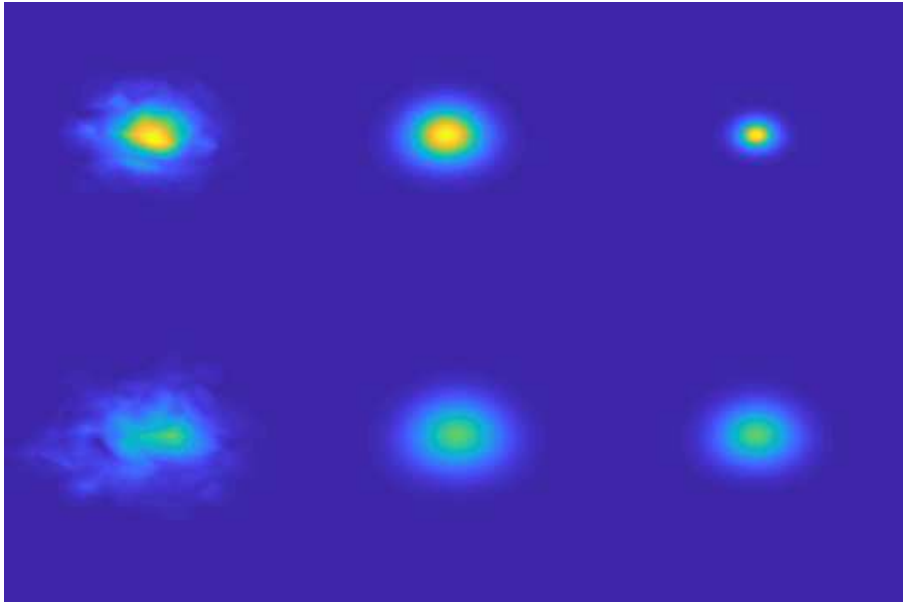}}
\end{picture}
\caption{\label{fig:ampreconExample1} Simulated FFT amplitude, $f=0.2Hz$, of sources A2 and B1 on the left. Best fit using the amplitude information in the middle. Best fit using the phase information on the right. The lower row corresponding to the source B1 which is situated more deep is multiplied by 10 to be visible in the same color map.}
\end{figure}

\begin{figure}
\begin{picture}(120,100)
\epsfxsize=12cm
\put(0,0) {\epsffile{\imagepath 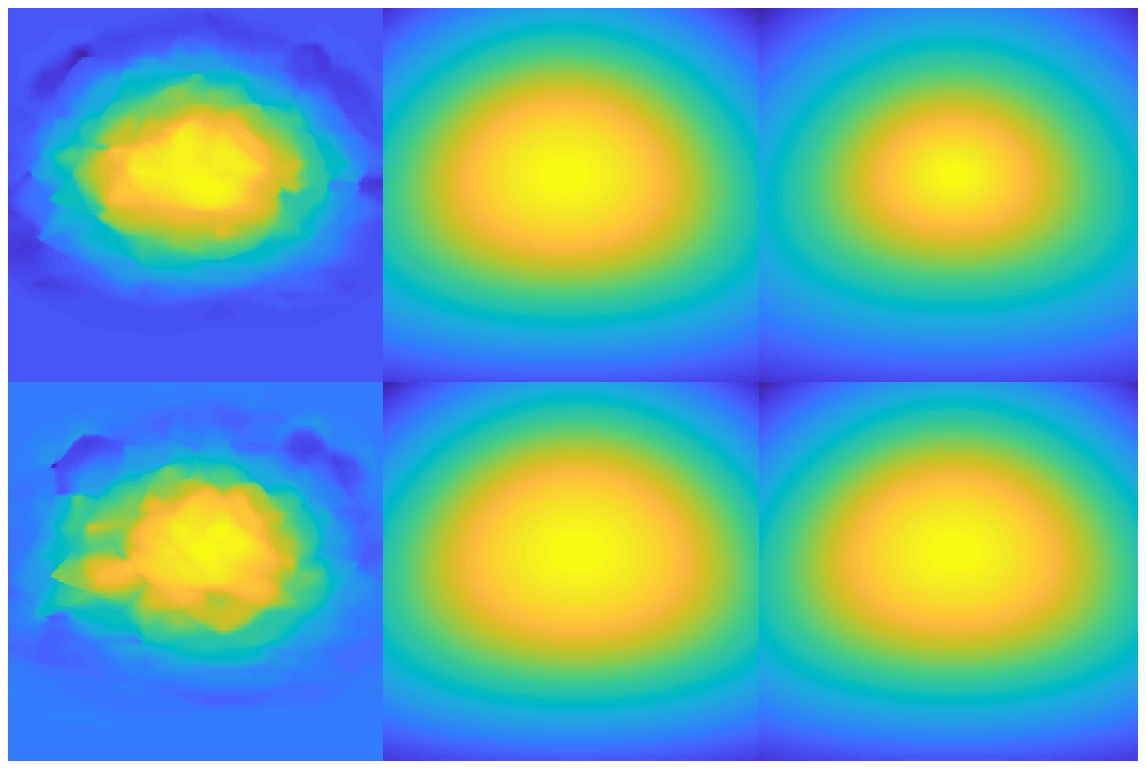}}
\end{picture}
\caption{\label{fig:phasereconExample1} Simulated FFT phase, $f=0.2Hz$, of sources A2 and B1 on the left. Best fit using the amplitude information in the middle. Best fit using the phase information on the right. All of these phase images have been normalized to show their shape in relation to one another.}
\end{figure}

In tables \ref{results2A} and \ref{results2B} are the results of the dynamic simulations. The location columns have the distance from the reconstructed source center to the true source center in same units as the simulations, both using the amplitude and phase information. The spectrum columns give the reconstructed source modulation amplitude compared to the non-noisy reconstruction of source A1, which again should be the easiest to reconstruct. In figures \ref{fig:ampreconExample1} and \ref{fig:phasereconExample1} we see examples of how close the model computed FFT amplitudes are to the simulated FFT amplitudes. What is not presented here is the full reconstructed spectra; the cutoff parameter $At$ was efficiently chosen so that all other frequency components other than the third FFT bin corresponding to exactly $f=0.2Hz$ were zero, in all reconstructions.

\begin{table}
\begin{tabular}{c|c|c}
          &    location error      &  spectrum magnitude S \\
center A  &   (amp) $\quad$ (phase)& (amp) $\quad$ (phase) \\
\hline
no noise(A1)   &   6.1  $\quad$  5.0     &      2.2  $\quad$  2.2 \\
no noise(A2)   &   6.1  $\quad$  5.0     &      1.8  $\quad$  1.8 \\
noise 1\% (A1) &   6.1  $\quad$  8.2     &      2.2 $\quad$   4.0 \\
noise 3\% (A2) &   6.1  $\quad$  8.2     &      1.8  $\quad$  3.2
\end{tabular}
\caption{\label{results2A} Dynamic reconstructions for the source location A, two modulation levels 2.2 and 1.8 and two simulated noise levels.}
\end{table}

\begin{table}
\begin{tabular}{c|c|c}
          &    location error      &  spectrum magnitude S \\
center B  &   (amp) $\quad$ (phase)& (amp) $\quad$ (phase) \\
\hline
no noise(B1)   &   3.8  $\quad$  1.5     &      5.2   $\quad$  438 \\
no noise(B2)   &   3.8  $\quad$  2.3     &      4.2   $\quad$  2409 \\
noise 1\% (B1) &   5.3  $\quad$  6.8     &      13.4  $\quad$   411 \\
noise 3\% (B2) &   5.3  $\quad$  8.1     &      11.0  $\quad$  14945
\end{tabular}
\caption{\label{results2B} Dynamic reconstructions for the source location B, two modulation levels 2.2 and 1.8 and two simulated noise levels.}
\end{table}

\section{Conclusions}\label{sec:conclusions}
In this article we have presented a new approach to detecting static and time-varying point-like heat sources. The numerical work here can be thought of as a proof of concept of this approach, even if there is a lot of work to do in order to develop the algorithms, especially as the data changes from simple and simulated to realistic and noisy. This simple model can in principle be used to approximate more complicated cases.

The numerical simulations turned out to be very challenging as we would need high numerical accuracy in 4D. In the future we hope to revisit the simulations with realistic parameters. As can be seen in figures \ref{fig:ampreconExample1} and \ref{fig:phasereconExample1} on the left, the simulations are not looking smooth - this could be improved by using a better and finer FEM mesh and possibly by adjusting FEM solver parameters. This is a challenging and laborous task in itself.

On the other hand, the example model fits of figures \ref{fig:reconExample},\ref{fig:ampreconExample1} and \ref{fig:phasereconExample1} are good enough to label this approach as promising. The simulation was done on a finite half ball with a semi-realistic boundary condition used in real bio-heat models; still the analytical point-source model using a solution in the infinite 3D space fits reasonably well. Another future project would be to further investigate ways to either modify our model to accommodate realistic features such as, for example, blood convection, or to analytically compute the different errors it induces when used in realistic scenarios. For example, the boundary condition \eqref{boundaryCondition} could perhaps be addressed analytically as well as different source shapes. In fact, any source could be thought of as a collection of point sources, each contributing as a simple model solution of \eqref{modelSolution1} and \eqref{modelSolution2}.

One very interesting future research subject is the efficient use of phase information. As the model \eqref{modelSolution1} predicts, the phase profile will be not attenuated exponentially, but linearly, although the phase is related to the amplitude and we will not be able to compute it more accurately from the data.

The results of tables \ref{results1}, \ref{results2A} and \ref{results2B} have to be taken critically because of the accuracy problems of the simulations and inherent complexity of the algorithms. Still, one could argue that the basic trends of an inverse problem are present: the more deeply situated source, the more difficult it is to reconstruct, and the more noise, the results get worse quickly. The dynamic results might look perplexing, but we think that the following explains them:
\begin{itemize}
\item Using the amplitude information, the location detection is almost unchanged as noise is increased, in all of the simulations. What is not shown here is the results for higher noise levels, from around 5\% onwards the reconstruction break down very suddenly. Also the noise presents itself differently as we compute FFT from the noisy data. This is another issue completely and is not covered here, but in general it is possible to detect very small periodic phenomena from noisy data if the sampling time is large. In our simulation the sampling time was quite average ten seconds. 
\item For some of the simulations, especially using the phase information, the reconstructed spectrum amplitude value change radically. These values are explained by a best fitting source location around 10 units towards or away from the true location, and so the source amplitude is then computed with or without this extra attenuating layer. This behaviour really highlights the ill-posedness of the problem, and to make the algorithms more usable we would need more robust handling of data.
\item The dynamic results of the more deeply located source B are better than of source A and even the static cases (based on the location error). This can be explained by the fact that the more deeply located source is actually more like a point source and thus the analytic model could work better in that case. Also, the dynamic cases use a lot more data for the source detection, and so if the FFT can be computed properly like in our simulations, it makes sense that the location is computed more accurately.
\item The reconstructions using the phase information might suffer from the more difficult algorithm implementation and also from larger deviation between the simulation and the analytical point source model, as the phase profile is not so attenuated.
\end{itemize}

\end{document}